\documentclass[]{article}

\usepackage[utf8]{inputenc}
\usepackage[T1]{fontenc}
\usepackage[english]{babel}
\usepackage{setspace}
\usepackage{pstricks}
\usepackage{pst-node}
\usepackage{amsmath}
\usepackage{amssymb}
\usepackage{graphicx}
\usepackage{fullpage}
\usepackage[bookmarks=true]{hyperref} 
\usepackage{tikz}
\usepackage{mathtools}
\usepackage[text={15cm,24cm}]{geometry}
\usepackage{enumerate}
\usepackage{caption}
\usepackage{subcaption}
\usepackage{algorithmicx}
\usepackage{algpseudocode}
\usepackage{algorithm}
\usepackage{color}
\usepackage{stmaryrd}

{~\ \hfill $\Box$\vspace{0,5cm}}

{~\ \hfill$\Diamond$\vspace{0,5cm}}

\newtheorem{prop}{Property}[section]

\newtheorem{rmk}{Remark}[section]

\begin{document}
	\title{Adaptive Network Flow with $k$-Arc Destruction}
	
	\author{
		T.\ Ridremont \footnotemark[1]
		\and
		D.\ Watel \footnotemark[2]
		\and
		P.-L. Poirion \footnotemark[3]
		\and
		C.\ Picouleau \footnotemark[4]
	}
	\date{}
	
	\def\thefootnote{\fnsymbol{footnote}}
	
	\footnotetext[1]{ \noindent  Conservatoire National des Arts et M\'etiers, CEDRIC laboratory, Paris (France). Email: {\tt
			thomas.ridremont@cnam.fr}
	}
	
	\footnotetext[2]{ \noindent  ENSIIE, SAMOVAR, Evry (France). Email: {\tt
			dimitri.watel@ensiie.fr}
	}
	
	\footnotetext[3]{ \noindent
		mathematical and algorithmic sciences lab, france research center, huawei technologies. Email: {\tt
			kiwisensei@gmail.com}
	}
	
	\footnotetext[4]{ \noindent
		Conservatoire National des Arts et M\'etiers, CEDRIC laboratory, Paris (France). Email: {\tt
			christophe.picouleau@cnam.fr}
	}

	\maketitle 
	
	\begin{abstract}
		When a flow is not allowed to be reoriented the Maximum Residual Flow Problem with $k$-Arc Destruction is known to be $NP$-hard for $k=2$. We show that when a flow is allowed to be adaptive the problem becomes polynomial for every fixed $k$.  
		
		\vspace{0.2cm}
		\noindent{\textbf{Keywords}\/}: adaptive flow, maximum flow, linear program, complexity.
	\end{abstract}

	\parindent=0cm

	\section{Introduction and motivation}
	The Maximum Residual Flow Problem with Arc Destruction is posed as follows. Given a directed network $G$ and a positive integer $k$ find a maximum flow such that when $k$ arcs are deleted from $G$ the value of a maximum flow in the residual network is maximum. Aneja {\it et al.} \cite{Aneja} and Du {\it et al.} \cite{Du} studied this problem when the flow cannot be reoriented. In \cite{Aneja} it is proved that the problem is polynomial time solvable for $k=1$ whereas  it is proved to be $NP$-hard for $k=2$ in  \cite{Du}.
	
	The Adaptive Maximum Flow Problem was first adressed by Bertsimas {\it et al.} in \cite{Bertsimas} as an alternative to the robust maximum flow problem. In this model, the flow can be adjusted after the arc failures occurred: the residual flow is the maximum flow in the network where the destroyed arcs are removed and where the capacity in each remaining arc equals the original flow of that arc. As described in \cite{Bertsimas} this adaptability is motivated by the pipeline applications where it is conceivable in practice to reroute the flow after the arc failures have occurred. In this note we are concerned with the Adaptative Maximum Residual Flow Problem, that is the maximum Residual Flow Problem with Arc Destruction when the flow can be reoriented. This problem is shown to be $NP$-hard when the number of arc destructions is a part of the instance \cite{Bertsimas}. We show that contrary to the case where the flow cannot be reoriented the Adaptative Maximum Residual Flow with $k$-Arc Destruction Problem is polynomial time solvable for any fixed $k$.
	
	The paper is organized as follows: the notations and the definition of the Adaptive Maximum Residual Flow with $k$-Arc Destruction Problem  (\emph{$k$-AMRFP}) are given in the next section. In section \ref{PL} we prove that \emph{$k$-AMRFP} is polynomial when the number of arc destructions is fixed. This is done by using a linear programming formulation.  We conclude and give some open questions in Section \ref{conclusion}.

	\section{Notations and Definitions}\label{notations}
	Let $G=(V,A,c)$ be a directed network with nonnegative capacities on arcs $c\ :\ A\rightarrow R$, $\vert V\vert=n$ vertices and $\vert A\vert=m$ arcs. Let $M$ be the $n\times m$ incidence matrix of $G$. Let $s,t\in V$ be the source and the sink of $G$, respectively. For convenience $A$ contains the dummy arc $(t,s)$ with unbounded capacity $c_{ts}=\infty$. Let $\varphi \ :\ A\rightarrow R$ be a (feasible) flow from $s$ to $t$ that satisfies both the capacity constraints and the flow conservation: $0\le\varphi\le c$ and $M\cdot\varphi=0$. Let $\varphi_{ts}$ be the value of $\varphi$. We say that the flow $\varphi$ is maximum if its value $\varphi_{ts}$ is maximal. The reader is referred to \cite{Ahuja} for all other aspects concerning flows in networks. 
	
	Let $A'\subset A\setminus\{(t,s)\}$ be a subset of failing arcs. Let $G'=(V,A\setminus A',c)$ be the residual network obtained from $G$ by deleting the failing arcs. Let $\varphi$ be a maximum flow for $G$. Then $\varphi'$ is an adaptive flow for $G'$ with respect to $\varphi$ if for each arc of $A\setminus A'$ the quantity over it is no more than in $\varphi$:  $M'\cdot\varphi'=0$ and $0\le\varphi'\le\varphi\le c$. 
	
	Let $k$ be a positive integer. We consider that the number of destroyed arcs is fixed to $k$, and that $A'$ can be any subset of $k$ arcs. Then the problem \emph{$k$-AMRFP}  is defined as follows:\\
	
	\textbf{Adaptive Maximum Residual Flow with $k$-Arc Destruction Problem (\emph{$k$-AMRFP})}
	
	\textit{INSTANCE: } A  directed network $G=(V,A,c)$ with a source $s \in V$ and a sink $t \in V$.    
	
	\textit{OBJECTIVE: } Find a maximum flow $\varphi$ on $G=(V,A,c)$ minimizing the loss of flow over all residual networks, that is $\min_{\varphi}\max\limits_{A'\subset A,\vert A'\vert=k}\{\varphi_{ts}-\varphi_{ts}'\}$.\\
	
	\begin{rmk}
		As shown by Aneja {\it et al.} in \cite{Aneja} (Lemma 2) if the initial flow on $G$ is not maximum then the maximum flow on the residual network cannot be better than the one obtaining when starting with a maximum flow on $G$. This the reason why we are only  searching for a maximum flow on $G$. 
	\end{rmk}
	
	The \emph{$k$-AMRFP} can be formulated as a two-person game, where a defender computes a feasible flow $\varphi$, the attacker then deletes $k$ edges, the defender finally computes $\varphi'$ a maximum adaptive flow in the remaining network. The objective of the defender is $\min_{\varphi} \max\limits_{A'\subset A,\vert A'\vert=k}\{\varphi_{ts}-\varphi_{ts}'\}$.\\

	\section{A polynomial size linear program for $k$-Arc Destruction}\label{PL}
	We show how \emph{$k$-AMRFP} can be formulated as a linear program. Suppose the arcs of $G$ are labeled $a_1,\ldots, a_m,a_{m+1}$, with  $a_{m+1}$ the dummy arc $(t,s)$. Thus in $M$ the $n\times (m+1)$ incidence matrix of the network $G$ the column $i$ corresponds to $a_i$. Let ${\cal A}=\{A_1,\ldots,A_{m\choose k}\}$ be the set consisting of the $m\choose k$ $k$-subsets  of $A\setminus\{a_{m+1}\}$. For each $i, 1\le i\le {m\choose k}$, we denote by $M^i$ the submatrix obtained from $M$ by deleting the $k$ columns corresponding to $A_i$ and by $\varphi^i\in R^{m+1-k}$ a $(m+1-k)$-dimensional vector.
	
	It is straightforward to verify that \emph{$k$-AMRFP} is modeled as the following integer program.
	
	\[
	\begin{aligned}
	\max\quad\varphi_{ts}+\theta   \\
	M\cdot\varphi & \quad = \quad 0 &\\
	\varphi & \quad \leq \quad c &\\
	M^i\cdot\varphi^i & \quad = \quad 0 &\quad 1\le i\le {m\choose k}\\
	\varphi^i & \quad \leq \quad \varphi &\quad 1\le i\le {m\choose k}\\
	\varphi^i_{ts} & \quad \geq \quad \theta &\quad 1\le i\le {m\choose k}\\
	\theta,\varphi,\varphi^i &\quad \geq \quad  0 \qquad & 1\le i\le {m\choose k}
	\end{aligned}
	\]
	
	Indeed by the $m\choose k$ constraints $M^i.\varphi^i=0$ we consider all adaptive flows resulting of  the destruction of any $k$ arcs. Moreover combining the constraints $\varphi^i_{ts}  \geq \theta$ with the objective function, the minimum value among these $m\choose k$ flows is maximal. Note also that the number of variables and constraints is $O(m^{k+1})$. Thus we have the following. 
	\begin{prop}
		\emph{$k$-AMRFP} is polynomial time solvable for any fixed integer $k>0$.
	\end{prop}

	\section{Conclusion and discussion}\label{conclusion}
	We proved that contrary to the non reoriented case the adaptive residual flow problem is polynomial when $k$ the number of destructed arcs is fixed. This is done using linear programming.

	A challenging problem is to find a polynomial time combinatorial algorithm to solve this problem. A natural way   to reach this goal is to try to adapt one of the classical maximum flow algorithms (see \cite{Ahuja} for such algorithms). The following example tries to explain why this task does not seem  to be easy. 
	
	Let $G$ be the network with three vertices $s,v,t$ with $b$ arcs of capacity $1$ from the source $s$ to $v$ and two arcs of capacity $b$ between $v$ to the sink $t$.  Consider the problem with $k=1$ arc destruction. The flow $\varphi$ where there is $1$ unit of flow for all the $(s,v)$ arcs and $b/2$ units of flow for the two $(v,t)$  arcs is  an optimal solution of \emph{$k$-AMRFP}. Since $k=1$ the problem consists in finding a maximum flow where the maximal value of a flow through an arc is minimal (see \cite{Aneja}).  This flow can be found by the algorithm given in \cite{Aneja} which uses a Newton's method. It seems that  this algorithm cannot be generalized  for $k\geq 2$. A second difficulty is the following:  most of the  standard combinatorial algorithms for maximum flows based on the structure of $G$ give either a flow $\phi_1$ where one of the two $(v,t)$ arcs has a flow $b$ and the second has a flow null or a flow $\phi_2$ with the opposite flows in the two $(v,t)$ arcs.  In our example the optimal solution $\varphi$ consists of a linear convex combination of both $\phi_1$ and $\phi_2$: $\varphi={1\over 2}(\phi_1+\phi_2)$. Even in this small example it seems not easy to design a combinatorial algorithm that  finds a maximum flow $\varphi$ minimizing $\max\limits_{a \in A}(\varphi(a))$.
	
	Other open problems concern several complexity aspects of  the Adaptive Maximum Flow Problem: as noted by Bertsimas {\it et al.} \cite{Bertsimas} determining $\varphi_{ts}'$, the maximum value of a residual flow, is $NP$-hard. However building an optimal initial flow (among the maximum flows) may be an easier problem if we do not include the value of the associated residual flow in the solution. The same kind of difficulties occur concerning the parameterized complexity.


\begin{thebibliography}{99}
		
		\bibitem{Ahuja}
		R.K. Ahuja, T. L.  Magnanti, J.B. Orlin, \textit{Networks flows: Theory, Algorithm, and Applications}, Prentice Hall (1993).
		
		\bibitem{Aneja}
		Y.P. Aneja, R. Chandrasekaran, K. P. K. Nair, \textit{Maximizing Residual Flow Under an Arc Destruction}, Networks, 38(4)  (2001) 194-198.
		
		\bibitem{Bertsimas}
		D. Bertsimas, E. Nasrabadi, S. Stiller, \textit{Robust and Adaptive Network Flows}, Operations Research, 61(5) (2013) 1218-1242. 
		
		\bibitem{Du}
		D. Du, R. Chandrasekaran, \textit{The Maximum Residual Flow Problem: NP-hardness with Two-Arc Destruction}, Networks, 50(3)(2007) 181-182. 
		
		\bibitem{Guo}
		J. Guo, Y. R. Shrestha, \textit{Parameterized Complexity of Edge Interdiction Problems}, COCOON2014, LNCS 8591 (2014) 166-178. 
		
	\end{thebibliography}
\end{document}